\documentclass[12pt]{article}

\usepackage[reqno]{amsmath}
\usepackage{amssymb, amsthm}
\usepackage{enumerate,arydshln} 
\usepackage{mathrsfs}
\newtheoremstyle{plainsl}%
	{\topsep}
	{\topsep}
	{\slshape} 
	{}
	{\normalfont\bfseries}
	{.}
	{ }
	{}

\theoremstyle{plainsl}
\newtheorem{theorem}{Theorem}[section]
\newtheorem{lemma}[theorem]{Lemma}
\newtheorem{corollary}[theorem]{Corollary}
\newtheorem{proposition}[theorem]{Proposition}

\DeclareMathOperator\fix{Fix}
\DeclareMathOperator{\BigO}{O}

\renewcommand\proof{\noindent\textsl{Proof. }}
\newcommand\sqr[2]{{\vbox{\hrule height.#2pt
    \hbox{\vrule width.#2pt height#1pt \kern#1pt
        \vrule width.#2pt}\hrule height.#2pt}}}

\renewcommand\qed{%
	\ifmmode\eqno\sqr53
	\else\nolinebreak\ \hfill\sqr53\medbreak\fi}

\numberwithin{equation}{section}

\newcommand{\erd}{Erd\H{o}s}

\newcommand\sym{\mathrm{Sym}(n)}
\newcommand\fip{{\mathcal A}}
\newcommand\ijfix{\triangleleft_{ij}}
\newcommand\comp{\mathcal C_{i,j}}
\newcommand\pr{\,^\prime}
\newcommand\up{\mathscr U_{p}}
\newcommand\rpart[1]{\mathscr R_{#1}}
\newcommand\rprime[1]{\mathscr R\,^\prime_{#1}}

\title{The exact bound for the Erd\H{o}s-Ko-Rado theorem for $t$-cycle-intersecting permutations}
\author{Karen Meagher \footnote{Research supported by NSERC Discovery Grant 341214-08.}\\
\small  Department of Mathematics and Statistics \\[-0.8ex]
\small University of Regina,  Regina, Saskatchewan, Canada\\[-0.8ex]
\small \texttt{kmeagher@math.uregina.ca}\\[+0.8ex]
Alison Purdy \footnote{Research supported by NSERC Postgraduate Scholarship.}\\
\small  Department of Mathematics and Statistics \\[-0.8ex]
\small University of Regina,  Regina, Saskatchewan, Canada\\[-0.8ex]
\small \texttt{purdyali@math.uregina.ca}
}

\begin{document}
\maketitle

\abstract{In this paper we adapt techniques used by Ahlswede and Khachatrian in their proof of the Complete \erd-Ko-Rado Theorem to show that if $n \geq 2t+1$, then any pairwise $t$-cycle-intersecting family of permutations  has cardinality less than or equal to $(n-t)!$.  Furthermore, the only families attaining this size are the stabilizers of $t$ points, that is, families consisting of all permutations having $t$ 1-cycles in common.  This is a strengthening of a previous result of Ku and Renshaw who proved the maximum size and structure of a $t$-cycle-intersecting family of permutations  provided $n \geq n_{0}(t)$ where $n_{0}(t)=O(t^2)$ and supports a recent conjecture by Ellis, Friedgut and Pilpel concerning the corresponding bound for $t$-intersecting families of permutations.}

\section{Introduction}
The \erd-Ko-Rado theorem~\cite{MR0140419} has been the focus of a great deal of interest since it was first published in 1961.  This theorem describes the size and structure of the largest collection of $k$-subsets of an $n$-set having the property that any two subsets have at least $t$ elements in common.  Such a collection is called \textsl{$t$-intersecting}.  One statement of the theorem is as follows:

\begin{theorem}\label{thm:fullekr}
  Let $t\leq k \leq n$ be positive integers.  Let $\mathcal F$ be a family of
  pairwise $t$-intersecting $k$-subsets of $\{1,\dots,n\}$.  There
  exists a function $n_{0}(k,t)$ such that for $n \geq n_{0}(k,t)$,
  \[|\mathcal F| \leq \binom{n-t}{k-t}.\] Moreover, for $n > n_{0}(k,t)$, 
  $\mathcal F$ meets this bound if and only if $\mathcal F$ is the collection
  of all $k$-subsets that contain a fixed $t$-subset.
\end{theorem}

In their original paper, \erd, Ko and Rado showed that $n_{0}(k,1)=2k$ and gave $t+(k-t)\binom{k}{t}^3$ as a necessary, although not optimal, lower bound on $n$ for values of $t$ other than 1.  In 1976, Frankl~\cite{MR519277} proved that $n_{0}(k,t)=(t+1)(k-t+1)$ for $t \geq 15$ and in 1984, Wilson~\cite{MR771733} showed that this bound was applicable for all values of $t$.  In 1997, Ahlswede and Khachatrian~\cite{MR1429238} published what is known as the Complete \erd-Ko-Rado Theorem which gives the maximum size and structure of $t$-intersecting set systems for all values of $n$, $k$ and $t$ and which includes a different proof of the exact lower bound on $n$. 

Let $\sym$ denote the symmetric group on $[n]=\{1,\dots,n\}$.  A family of permutations, $\fip \in \sym$, is \textsl{$t$-intersecting} if for any two permutations $\sigma,\pi \in \fip$, there exists a set $B=\{b_{1},\dots,b_{t}\} \subset [n]$ such that $\sigma(b_{i})=\pi (b_{i})$ for $i=1,\dots,t$.  In 1977, Deza and Frankl~\cite{MR0439648} proved that a $1$-intersecting family will have size at most $(n-1)!$.  Significantly later, a number of different approaches were used to show that the only families achieving this size were the cosets of stabilizers of a point (see~\cite{MR2009400},~\cite{MR2489272},~\cite{MR2061391},~\cite{MR2419214}).  

More recently, Ellis, Friedgut and Pilpel~\cite{MR2784326} made a significant advance by proving a long standing conjecture from~\cite{MR0439648} regarding the maximum size of a $t$-intersecting family of permutations with the following theorem.
\begin{theorem}\label{ellis}
For any given $t \in \mathbb N$ and $n$ sufficiently large relative to $t$, if $\fip$ is a $t$-intersecting family of permutations from $\sym$, then $\left | \fip \right| \leq (n-t)!$, with equality if and only if there exist $t$ distinct integers, $a_{1},\dots,a_{t} \in [n]$, , and $t$ distinct integers, $b_{1},\dots,b_{t} \in [n]$,  such that $\sigma(a_{i})=b_{i}$ for all $\sigma \in \fip$ and $i=1,\dots,t$.
\end{theorem}  
Their proof used eigenvalue methods and the representation theory of $\sym$.  In this paper we adapt the combinatorial method of Ahlswede and Khachatrian to prove a similar result for $t$-cycle-intersecting permutations.  

The concept of $t$-cycle-intersection was introduced in an earlier paper by Ku and Renshaw~\cite{MR2423345}.  A family of permutations is \textsl{$t$-cycle-intersecting} if any two permutations in $\fip$, when written in their cycle decomposition form, have at least $t$ cycles in common.  If $\fip$ is $t$-cycle-intersecting, it is $t$-intersecting but the converse is not true.  They proved the following theorem. 
\begin{theorem}\label{KR}
Suppose $\fip \subseteq \sym$ is $t$-cycle-intersecting and $n \geq n_{0}(t)$ where $n_{0}(t)=\BigO(t^{2})$.  Then $\left| \fip \right| \leq (n-t)!$ with equality if and only if $\fip$ is the stabilizer of $t$ points.
\end{theorem}

Our main result is a refinement of Theorem~\ref{KR}.

\begin{theorem}\label{KR with bound}
Suppose $\fip \subseteq \sym$ is $t$-cycle-intersecting.  If $n \geq 2t+1$, then $\left| \fip \right| \leq (n-t)!$ with equality if and only if $\fip$ is the stabilizer of $t$ points.
\end{theorem}

Although the size and structure of the largest $t$-cycle-intersecting family of permutations is implied by Theorem~\ref{ellis}, our proof gives the exact lower bound on $n$.  This bound, $n \geq 2t+1$, is the same as that implied by a conjecture of Ellis, Friedgut and Pilpel in~\cite{MR2784326} for $t$-intersecting families of permutations; their conjecture remains to be proved and they state the opinion that the proof will require new techniques.  However, the method used in this paper has been used in~\cite{alison} to prove that this bound holds when the $t$-intersecting families of permutations are restricted to those closed under the fixing operation introduced by Cameron and Ku in~\cite{MR2009400}.

Determining the exact lower bound on $n$ is an important step towards determining the largest $t$-cycle-intersecting families of permutations for all values of $n$ and $t$.  To show that our bound is the best possible, in Section 4 we present a $t$-cycle-intersecting family of permutations that is larger than the stabilizer of $t$ points for all $n$ such that $t+3 < n < 2t+1$.

\section{Preliminary results}


In their proof of the Complete \erd-Ko-Rado Theorem, Ahlswede and Khachatrian make use of properties of set systems that are closed under the well known left-shifting operation.  (For information about this operation, see~\cite{MR905277}.)  In Sections 2.1 and 2.2, we consider two operations on families of permutations. The first,  $ij$-fixing, increases the number of elements fixed by the permutations.  The second operation, compression, acts on the fixed points of the permutations in a manner similar to the left-shifting operation for set systems.  

The proof of Lemma~\ref{main lemma} relies on properties of families of permutations that are closed under these two operations.  The results presented in Sections 2.1 and 2.2 (except for Proposition~\ref{KR fix}) are needed to extend this result to all $t$-cycle-intersecting families of permutations.  They show that any maximal $t$-cycle-intersecting family can be transformed into a family that is closed under both operations by performing $ij$-fixing operations followed by compression operations and that the resulting family has the same size as the original family, is $t$-cycle-intersecting and is the stabilizer of $t$-points only if the original family was the stabilizer of $t$-points.  

The propositions in Section 2.2 generally follow from the definitions of the two operations although the proofs often require a somewhat tedious examination of multiple cases.  For convenience, we will use $I(n,t)$ to denote the collection of all $t$-cycle-intersecting families of permutations.  

\subsection{Fixing}

This operation was introduced by Ku and Renshaw in~\cite{MR2423345} and we use their notation.

For a permutation, $\sigma$, the \textsl{$ij$-fixing}, $_{[ij]}\sigma$, is defined as follows:

\begin {enumerate}
\item  if $\sigma(i) \neq j$, then $_{[ij]}\sigma=\sigma$; 
\item  if $\sigma(i)= j$, then $_{[ij]}\sigma(x)=
\begin {cases}
	i & \mbox {if } x=i \mbox { ,} \\
	j & \mbox {if } x=\sigma^{-1}(i) \mbox { ,} \\
	\sigma(x) & \mbox {otherwise} \mbox {.} \\
	\end {cases} $
\end {enumerate}
Further, for a set of permutations, $\fip$, the $ij$-fixing of $\fip$, $\ijfix(\fip)$, is defined as

\begin{itemize}
\item[] $\triangleleft_{ij}(\fip)=\{\triangleleft_{ij}(\sigma)\: : \: \sigma \in \mathcal A\}$, where
\item[] $\triangleleft_{ij}(\sigma)=\begin{cases} _{[ij]}\sigma & \mbox { if } _{[ij]}\sigma \notin \mathcal A \mbox {,} \\
\sigma & \mbox { otherwise}\mbox {.}\\
\end{cases}$
\end{itemize}



We will refer to a family of permutations, $\fip$, having the property that $\triangleleft_{ij}(\fip)=\fip$ for all $i \neq j \in [n]$ as a \textsl{fixed} family.  Ku and Renshaw show in~\cite{MR2423345} that any $t$-cycle-intersecting family of permutations can be transformed into a fixed $t$-cycle-intersecting family of permutations by repeated applications of the $ij$-fixing operation and that the number of distinct permutations in the fixed family will be the same as the number in the original family.  They also prove the following propositions.   

\begin{proposition}\label{KR stab} 
Let $n > t+1$ and let $\fip \in I(n,t)$.  If $\ijfix(\fip)$ is the stabilizer of $t$ points, then so is $\fip$.
\end{proposition}

For a permutation, $\sigma$, let $\fix(\sigma)=\{x: \sigma(x)=x\}$.  For a family of permutations, $\fip$, let $\fix(\fip)=\{\fix(\sigma): \sigma \in \fip\}$.  

\begin{proposition}\label{KR fix}
Suppose $\fip \in I(n,t)$ is fixed.  Then $\fix(\fip)$ is a $t$-intersecting family of subsets of $[n]$.

\end{proposition}


\subsection{Compression}
We now define an operation on families of permutations that is designed to make the fixed points in the permutations as small as possible while retaining the size and structure of the family.  We will use the term compression to describe this operation.

For a permutation, $\sigma \in \sym$, and $i,j \in [n]$ with $i<j$, the \textsl{(i,j)-compression} of $\sigma$, denoted by $\sigma_{i,j}$, is the permutation defined as:

\begin {enumerate}
\item  if $\sigma(i)=i$ or $\sigma(j) \neq j$, then $\sigma_{i,j}=\sigma$; 
\item  if $\sigma(i)\neq i$ and $\sigma(j)=j$, then $\sigma_{i,j}(y)=
\begin {cases} 
	i & \mbox {if } y=i \mbox { ,} \\
	\sigma(i) & \mbox {if } y=j \mbox { ,} \\
	j & \mbox {if } y=\sigma^{-1}(i) \mbox { ,} \\
	\sigma(y) & \mbox {otherwise}\mbox { .} \\
	\end {cases} $
\end {enumerate}

The compression operations does not change the size or number of cycles in the cycle decomposition of a permutation.  Its effect  on the set $\fix(\sigma)$ when $j$ is fixed under $\sigma$ and $i$ is not is to remove $j$ from $\fix(\sigma)$ and replace it with $i$.  Unlike the $ij$-fixing operation, the compression operation does not change the number of fixed points and will have the same effect on the set $\fix(\sigma)$ as applying the standard left shifting operation to $\fix(\sigma)$.

Let $\mathcal A$ be a collection of permutations.  We define the $(i,j)$-compression, $\comp(\mathcal A)$, as follows:
\begin{itemize}
\item[] $\comp(\mathcal A)=\{\comp(\sigma)\: : \: \sigma \in \mathcal A\}$, where
\item[] $\comp(\sigma)=\begin{cases} \sigma_{i,j}  & \mbox { if } \sigma_{i,j} \notin \mathcal A \mbox { ,}\\
\sigma  &\mbox { otherwise}\mbox { .}\\
\end{cases}$
\end{itemize}

We will refer to a family of permutations, $\mathcal A$, that has the property that $\comp(\mathcal A)=\mathcal A$ for all $1 \leq i < j \leq n$ as a \textsl{compressed} family.

A $t$-cycle-intersecting family of permutations, $\fip \subseteq \sym$, is \textsl{maximal} if for all $\sigma \in \sym \backslash \fip$, the family $\fip \cup \{\sigma\}$ is not $t$-cycle-intersecting.  The next four propositions show that any maximal fixed $t$-cycle-intersecting family of permutations can be transformed into a compressed fixed $t$-cycle-intersecting family with the same cardinality as the original family.  The final proposition in this section allows us to conclude that if the resulting compressed family is the stabilizer of $t$ points, then the original family was the stabilizer of $t$ points.       

\begin{proposition}\label{comp size}
Let $\mathcal A$ be any collection of permutations from $\sym$.  Then for any $i,j \in [n]$ with $i<j$, $$\left |\comp(\mathcal A) \right | = \left | \mathcal A \right |\:.$$  
\end{proposition}

\proof
It is clear from the definition of the compression operation that $\comp(\mathcal A)$ will not be larger than $\mathcal A$ and that the size will decrease only if there are two distinct permutations, $\sigma$ and $\pi$, in $\mathcal A$ such that $\sigma_{i,j}=\pi_{i,j}$ and $\sigma_{i,j} \notin \mathcal A$. 

Suppose that there are two such permutations.  Then $\sigma(i) \neq i$, $\pi(i) \neq i$, and $\sigma(j) = \pi(j) = j$.  Also, $\sigma_{i,j}(\sigma^{-1}(i))=j =\pi_{i,j}(\pi^{-1}(i))$ and thus it follows from $\sigma_{i,j}=\pi_{i,j}$ that $\sigma^{-1}(i)=\pi^{-1}(i)$.  From the definition of the compression operation we have that
$$\sigma(x) = \sigma_{i,j}(x)=\pi_{i,j}(x)=\pi(x)$$ for all $x \in [n] \backslash \{i, j, \sigma^{-1}(i)=\pi^{-1}(i)\}$. 

We have already seen that $\sigma(j)=\pi(j)$ and that $\sigma^{-1}(i)=\pi^{-1}(i)$ from which it follows that $\sigma(\sigma^{-1}(i))=\pi(\sigma^{-1}(i))$.  Since $\sigma$ and $\pi$ have been shown to agree on all elements of $[n]$ except $i$, they must agree on $i$ also.  Hence, if $\sigma_{i,j}=\pi_{i,j}$, then $\sigma=\pi$.
\qed

\begin{proposition}\label{comp}
Let $\mathcal A$ be any collection of permutations from $\sym$. If $\mathcal A$ is not compressed, then it can be transformed into a compressed family by the application of a finite number of compression operations.
\end{proposition}
\proof(For a more detailed treatment see Theorem 5.2.2 in~\cite{alison}.)\\
If $\mathcal A$ is not compressed, there exists some $\sigma \in \mathcal A$ and $i,j \in [n]$ with $i<j$  such that $\comp(\sigma) \neq \sigma$.  Applying the $(i,j)$-compression operation to $\mathcal A$ has the same effect on the set $\fix(\sigma)$ as applying the standard left shifting operation for sets to $\fix(\sigma)$.  If $\comp(\mathcal A)$ is not compressed, there exists some $\pi \in \comp(\mathcal A)$ and $i\pr<j\pr \in [n]$ and  such that $\mathcal C_{i\pr,j\pr}(\pi) \neq \pi$.  Applying the $(i\pr,j\pr)$-compression operation to $\comp(\mathcal A)$ will left shift the set $\fix(\pi)$.  The number of permutations in the family is finite, the number of fixed points in any permutation is finite and unchanged by the compression operation, and the fixed points are always shifted to smaller integers.  Thus this process is necessarily finite.  \qed

\begin{proposition}\label{comp int}
Let $\fip \in I(n,t)$ be fixed.  Then $\comp(\fip) \in I(n,t)$ for any $i,j \in [n]$ with $ i < j$.
\end{proposition}
\proof
Let $\sigma$ and $\pi$ be two distinct permutations in $\fip$.  Clearly, if either $\comp(\sigma)=\sigma$ and $\comp(\pi)=\pi$ or $\comp(\sigma) \neq \sigma$ and $\comp(\pi)\neq \pi$, 
then $\comp(\sigma)$ and $\comp(\pi)$ will be $t$-cycle-intersecting.  Therefore, assume that $\comp(\sigma)=\sigma$ and $\comp (\pi) \neq \pi$.  Then $\pi(j)=j$ and $\pi(i)\neq i$.  There are three possible reasons why $\sigma$ is not changed by the compression operation.  
\begin{enumerate}
\item[] Case 1:  $\sigma(i)=i$.\\
If $\sigma(j) \neq j$, then $\sigma$ and $\pi$ must have $t$ cycles in common which do not involve $i$ or $j$.  These $t$ cycles will also belong to $\comp(\pi)$.  If $\sigma(j)=j$ then $\sigma$ and $\pi$ must have at least $t-1$ common cycles not involving $i$ or $j$.  Then $\sigma$ and $\comp(\pi)$ will have these $t-1$ cycles and the cycle $(i)$ in common. 
\item[] Case 2:  $\sigma(j) \neq j$.\\
Let $\pi(i)=x$.  If $\sigma$ and $\pi$ have $t$ common cycles not involving $i$ or $j$, these cycles will also be in $\comp (\pi)$.  Clearly, $\sigma$ and $\pi$ cannot have a common cycle involving $j$.  If $\sigma$ and $\pi$ have a common cycle involving $i$, it follows $\sigma(i) \neq i$.  $\fip$ is a fixed family and $\pi(i) \neq i$, so the permutation $_{[ix]}\pi$ must be in $\fip$.  Since $_{[ix]}\pi(i)=i$ and $_{[ix]}\pi(j)=j$, it follows that $\sigma$ and $_{[ix]}\pi$ must have $t$ cycles in common not involving $i$ or $j$.  These $t$ cycles will also be in $\pi$ and $\comp(\pi)$.  
\item[] Case 3:  $\sigma(i) \neq i$, $\sigma(j) =j$ and $\sigma_{i,j} \in \fip$.\\
Since $\sigma_{i,j}(i)=i$ and $\sigma_{i,j}(j) \neq j$, the permutations $\sigma_{i,j}$ and $\pi$ must have $t$ cycles in common not involving $i$ or $j$.  These $t$ cycles will also belong to $\sigma$ and $\comp(\pi)$.
\qed   
\end{enumerate}  

\begin{proposition}\label{comp fixed}
Let $\fip \in I(n,t)$ be maximal and fixed.  Then for any $i,j \in [n]$, $i<j$, the family $\comp(\fip)$ is a fixed family of permutations.
\end{proposition}
\proof
Let $\sigma$ be an element of $\comp(\fip)$ for arbitrary $i<j \in [n]$.  It is sufficient to show that $_{[xy]}\sigma \in \comp(\fip)$ for any $x \neq y \in [n]$.  If $\sigma(x) \neq y$, then $_{[xy]}\sigma=\sigma \in \comp(\fip)$ so we will assume that $\sigma(x)=y$.  We will consider two cases each with several subcases.
\begin{enumerate}
\item[]Case 1:  $\sigma \notin \fip$.\\
This implies that there exists some permutation $\pi \in \fip$ such that $\pi_{i,j}=\sigma$.  Thus $\sigma(j)=\pi(i) \in [n]\backslash\{i,j\}$ and $\sigma(i)=i$ so $x \neq i$.  We will consider when $_{[xy]}\sigma$ fixes $j$ and when it does not fix $j$.
\begin{enumerate}
\item $_{[xy]}\sigma(j)=j$.\\
We claim that $_{[xy]}\sigma \in \fip$ and that $_{[xy]}\sigma(i)=i$ which means that $_{[xy]}\sigma  \in \comp(\fip)$. In order for $_{[xy]}\sigma(j)=j$, either $x=j$ or $(x j)$ is a $2$-cycle in $\sigma$.  If $x=j$, then $_{[xy]}\sigma=_{[iy]}\pi$ since $y=\sigma(j)=\pi(i)$.  If the $2$-cycle $(xj)$ is in $\sigma$, then $_{[xy]}\sigma=_{[jx]}\sigma=_{[ix]}\pi$.  In either case, it follows that $_{[xy]}\sigma  \in \comp(\fip)$ since both $_{[iy]}\pi$ and $_{[ix]}\pi$ fix $i$ and $\pi \in \fip$ which is a fixed family.  
\item $_{[xy]}\sigma(j) \neq j$.\\
This implies that the $2$-cycle $(xi)$ is not in $\pi$ and that $x \neq j$.  
Since $x \neq i$, it then follows from the definitions of the fixing and compression operations that $$(_{[xy]}\pi)_{ij}=_{[xy]}(\pi_{ij})=_{[xy]}\sigma$$ except when $y=j$.  If $y=j$, then $(_{[xi]}\pi)_{ij}=_{[xy]}\sigma$.  Again, since $\fip$ is a fixed family both $_{[xy]}\pi$ and $_{[xi]}\pi$ will be in $\fip$ and thus $_{[xy]}\sigma$ will be in $\comp(\fip)$.    

\end{enumerate}
\item[]  Case 2:  $\sigma \in \fip$.\\
If $_{[xy]}\sigma = (_{[xy]}\sigma)_{i,j}$, we're done since $_{[xy]}\sigma \in \comp(\fip)$ easily follows from the requirement that $\fip$ be a fixed family.  Therefore, we will assume that $_{[xy]}\sigma \neq (_{[xy]}\sigma)_{i,j}$ and so $_{[xy]}\sigma$ will be in $\comp(\fip)$ if and only if $(_{[xy]}\sigma)_{i,j}$ is in $\fip$. This assumption implies that $_{[xy]}\sigma(j)=j$ and that $_{[xy]}\sigma(i) \neq i$ and, in particular, that $x \neq i$.  We will consider when $\sigma \neq \sigma_{i,j}$ and when $\sigma=\sigma_{i,j}$ as separate cases.
\begin{enumerate}
\item $\sigma \neq \sigma_{i,j}$.\\
Since $\sigma$ is an element of both $\fip$ and $\comp(\fip)$, it follows that $\sigma_{i,j} \in \fip$.  Let $\sigma_{i,j}(x)=z$.  (Note that $z\!=\!y$ unless $y\!=\!i$; then $z\!=\!j$.)  Since $\fip$ is a fixed family, $_{[xz]}(\sigma_{i,j})$ will also be in $\fip$. 
But $_{[xz]}(\sigma_{i,j})=(_{[xy]}\sigma)_{i,j}$ and so $\comp(_{[xy]}\sigma)=_{[xy]}\sigma$ and $_{[xy]}\sigma \in \comp(\fip)$.

\item  $\sigma =\sigma_{i,j}$.\\
Since $_{[xy]}\sigma(i) \neq i$ implies that $\sigma(i)\neq i$, it follows from $\sigma=\sigma_{i,j}$ that $\sigma(j) \neq j$.  
Let $_{[xy]}\sigma(i)=z$ and consider the permutation $\pi=_{[iz]}(_{[xy]}\sigma)$.  Since $\pi, \sigma \in \fip$ and the cycles in $\sigma$ containing $i$ and $j$ will not be in $\pi$, there must be at least $t$ other cycles in $\sigma$ and these will also be in $\pi$.  Note that $(_{[xy]}\sigma)_{i,j}$ will also contain these cycles and thus will $t$-cycle-interesect with both $\sigma$ and $\pi$.  Let $\rho$ be any other permutation in $\fip$.  It must $t$-cycle-intersect with both $\sigma$ and $\pi$ and thus it will $t$-cycle-intersect with $(_{[xy]}\sigma)_{i,j}$.  Since $\fip$ is maximal, it follows that $(_{[xy]}\sigma)_{i,j}\! \in\! \fip$ and so $\comp(_{[xy]}\sigma)\!=\!_{[xy]}\sigma$.  \qed
\end{enumerate}
\end{enumerate}  
   
\begin{proposition}\label{comp stab}
Let $\fip \in I(n,t)$ and let $n > t+1$.  If $\comp(\fip)$ is the stabilizer of t points, then so is $\fip$.
\end{proposition}
\proof
Suppose $\fip$ is $t$-cycle-intersecting and $\comp(\fip)$ is the stabilizer of $t$ points.  If $\fip=\comp(\fip)$ the proposition is trivial so assume that $\fip \neq \comp(\fip)$.  Let $\{x_{1},x_{2},\dots,x_{n}\} = [n]$ and let $\comp(\fip)$ be the stabilizer of $x_{1},x_{2},\dots,x_{t}$.  Consider the permutation $(x_{t+1},\dots,x_{n})$.  Clearly, $\sigma \in \comp(\fip)$. 
\begin{itemize}
\item[]Case 1:  $ \sigma\in \fip$\\
Suppose there is a permutation $\pi \in \fip$ such that $(x_{t+1},\dots,x_{n})$ is a cycle in $\pi$.  Since $\pi$ must have $t$ cycles in common with $\sigma$, it follows that $\pi=\sigma$.  Hence $\sigma$ is the only permutation in $\fip$ containing $(x_{t+1},\dots,x_{n})$ and all other permutations must fix $x_{1},x_{2},\dots,x_{t}$.  From Proposition~\ref{comp size} we have that $\left |\fip  \right | = \left | \comp(\fip) \right|$, so $\fip$ must be the stabilizer of $t$ points.
\item[]Case 2: $ \sigma \notin \fip$\\
Since $\sigma \in \comp(\fip)$, there must be a permutation $\pi \in \fip$ such that $\pi_{i,j} = \sigma$.   Then $\pi =  (y_{t+1},\dots,y_{n})$ where $j \in \{y_{1},\dots, y_{t}\}$ and $i \in \{y_{t+1},\dots,y_{n}\}$ and $y_{k}=x_{k}$ for all $y_{k} \neq i,j$.  The argument from Case 1 can now be applied to show that $\fip$ is the stabilizer of the $t$ points $y_{1},\dots,y_{t}$.  \qed
\end{itemize}  

\subsection{Generating sets}

The concept of generating sets for set systems was introduced in~\cite{MR1429238}.  In order to define an analogous concept for permutations. we first define the \textsl{up-permutation} of a set, $B \subseteq [n]$, as $$\mathscr U_{p}(B) =\{\sigma \in \sym : B \subseteq \fix(\sigma)\}\:.$$  
For a collection of sets, $\mathcal B$, we define the up-permutation to be $$\mathscr U_{p}(\mathcal B) = \{\mathscr U_{p}(B):B \in \mathcal B \} \: .$$
Finally, a collection of sets, $g(\mathcal A)$, is a generating set for $\mathcal A \subseteq \sym$, if $g(\mathcal A)$ does not contain any sets of cardinality $n\!-\!1$ and $$\mathscr U_{p}(g(\mathcal A))= \mathcal A\;.$$  
The set of all generating sets of $\mathcal A$ is denoted by $G(\mathcal A)$.

Unlike $t$-intersecting set systems, where the set system itself is a generating set, not all $t$-cycle-intersecting families of permutations will have a generating set.  The following lemma proves the existence of generating sets for maximal fixed $t$-cycle-intersecting families.  

\begin{lemma}\label{existence}
Let $\fip \in I(n,t)$ be maximal and fixed.  Then $\fix(\fip)$ is a generating set for $\fip$.
\end{lemma}

\proof
By Proposition~\ref{KR fix}, the set system $\fix(\fip)$ is $t$-intersecting.  Hence $\up(\fix(\fip))$ is a $t$-cycle-intersecting family of permutations.  If $\sigma \in \fip$, then $\sigma \in \up(\fix(\fip))$ and so $\fip \subseteq \up(\fix(\fip))$.  Since $\fip$ is maximal, it follows that $\fip = \up(\fix(\fip))$ as required.
\qed

\begin{corollary}\label{not empty}
If $\fip \in I(n,t)$ is maximal and fixed, then $G(\fip) \neq \emptyset$.
\end{corollary}

\begin{lemma}\label{gen set int}
Let $\fip \in I(n,t)$ and let $n > t+1$.  Then any $g(\fip) \in G(\fip)$ is a $t$-intersecting set system.
\end{lemma}

\proof
Assume that $g(\fip) \in G(\fip)$ is not $t$-intersecting.  Then there exist $A,B \in g(\fip)$ such that $\left| A \cap B \right| < t $.  Since by the definition of generating sets neither $A$ nor $B$ can have cardinality $n-1$, there will be permutations in $\up(A)$ and $\up(B)$ that agree only on the elements in $A \cap B$ and hence are not $t$-cycle-intersecting.
\qed 

Some additional definitions are needed in order to describe the specific type of generating set required for the proof of the main theorem.

For $B = \{b_{1}, b_{2},\dots,b_{k}\} \in \binom{[n]}{k}$, where $b_{1}<b_{2}<\cdots<b_{k}$, the collection of all sets that can be obtained from $B$ by left-shifting is denoted by $\mathscr L(B)$ and is formally defined as follows:
$$\mathscr L(B) = \{A = \{a_{1}, a_{2},\dots,a_{k}\} \in \binom{[n]}{k} : a_{i} \leq b_{i} \mbox{ for all } i \in [k]\}\:.$$ 
For a $k$-set system, $\mathcal B \subseteq \binom{[n]}{k}$, $$\mathscr L(\mathcal B) = \{\mathscr L(B): B \in \mathcal B\}.$$  If $\mathscr L(\mathcal B) = \mathcal B$, then $\mathcal B$ is said to be \textsl {left-compressed}.  Unlike the usual left-shifting operation for collections of subsets, the number of sets in $\mathscr L(\mathcal B)$ may be greater than the number in $\mathcal B$.  However, a collection of sets will be left-compressed in this sense if and only if it is left-shifted in the usual sense. 

For a generating set $g(\mathcal A) \in G(\mathcal A)$, consider $\mathscr L(g(\mathcal A))$.  Now define $\mathscr L_{*}(g(\mathcal A))$ as the set of minimal (in the sense of set inclusion) elements of $\mathscr L(g(\mathcal A))$.  Let $G_{*}(\mathcal A)$ be the set of all generating sets of $\mathcal A$ such that $\mathscr L_{*}(g(\mathcal A))= g(\mathcal A)$.  By this definition, any generating set in $G_{*}(\mathcal A)$ will be left-compressed and minimal with respect to inclusion.

For a set, $B \subseteq [n]$, let $s^{+} (B)$ denote the largest element in $B$.  For a generating set, $g(\fip)$, let $s^{+}(g(\fip))=\max\{s^{+} (B) : B \in g(\fip)\}$.  Finally, let $s_{min}(G(\fip))=\min\{s^{+}(g(\fip)) : g(\fip) \in G(\fip)\}$.  

\begin{lemma}\label{smin geq t}

Let $\fip \in I(n,t)$ and let $n \geq t+1$.  Then $s_{min}(G(\fip)) \geq t$.  Furthermore, if $s_{min}(G(\fip))=t$, then $\fip$ is  the stabilizer of $t$ points.
\end{lemma}

\proof
Suppose that $s_{min}(G(\fip)) < t$.  Then there is some $g(\fip) \in G(\fip)$ such that $s^{+}(g(\fip)) < t$.  Let $B \in g(\fip)$. Clearly $\left| B \right | < t$ and so $\up(G(\fip))$ is not $t$-intersecting, a contradiction.  Now suppose that $s_{min}(G(\fip)) = t$.  Then there exists some $g(\fip) \in G(\fip)$ such that $s^{+}(g(\fip))=t$.  Since we have shown that $\left | B \right| \geq t$ for all $B \in g(\fip)$, it follows that $g(\fip)=\{[t]\}$ 
\qed    

The requirement  in the remaining lemmas in this section that $\fip \in I(n,t)$ be maximal, fixed and compressed is necessary to ensure that $G_{*}(\fip) \neq \emptyset$.
\begin{lemma}\label{properties}
Let $\fip \in I(n,t)$ be maximal, fixed and compressed and let $g(\fip)$ be any generating set of $\fip$.  Then:
\begin{enumerate}
\item $\mathscr L_{*}(g(\fip)) \in G(\fip)$;
\item $s^{+}(\mathscr L_{*}(g(\fip))) \leq s^{+}(g(\fip))$;
\item For any set, $B$, that can be obtained by left shifting a set in $g(\fip)$, either $B$ or a proper subset of $B$ will be in $\mathscr L_{*}(g(\fip))$.
\end{enumerate}
\end{lemma}
\proof
The last two statements follow easily from the definition of $\mathscr L_{*}(g(\fip))$.

To prove the first statement, take any $g(\fip) \in G(\fip)$.  Since $g(\fip) \subseteq \mathscr L(g(\fip))$, it follows from the definitions of the up-permutation operation and $\mathscr L_{*}(g(\fip))$ that $\fip \subseteq \up (\mathscr L_{*}(g(\fip)))$.

We now show that $\up (\mathscr L_{*}(g(\fip))) \subseteq \fip$. Let $\sigma \in \up(\mathscr L_{*}(g(\fip)))$.  Since $\up(\mathscr L_{*}(g(\fip))) \subset \up(\mathscr L(g(\fip)))$, there is some set $B \in \mathscr L(g(\fip))$ such that $B \subseteq \fix(\sigma)$.  If $B \in g(\fip)$, then $\sigma \in \fip$  so assume $B \notin g(\fip)$.  Then there is some $B\pr \in g(\fip)$ such that $B \in \mathscr L(B\pr)$.  Note that $\left| B\right|=\left| B\pr\right|$.

Let $B=\{x_{1},x_{2},\dots,x_{\ell}\}$ and let $B\pr=\{y_{1},y_{2},\dots,y_{\ell}\}$.  Then $x_{i} \leq y_{i}$ for all $i \in [\ell]$ and, since $B\notin g(\fip)$, there is at least one $i \in \ell$ such that $x_{i}\neq y_{i}$.

Let $\fix(\sigma)=\{x_{1},\dots,x_{\ell}\} \cup \{z_{1},\dots,z_{m}\}$ and choose a permutation, $\pi$, such that $\fix(\pi)=B\pr \cup Z$, where $Z=\{z_{1},\dots,z_{m}\}$.  Since $B\pr \in g(\fip)$, it follows that $\pi \in \fip$.  Note that $B\pr$ and $Z$ are not necessarily disjoint sets.  If $x_{1}= y_{1}$, then $\fix(\pi)=\{x_{1},y_{2},\dots,y_{\ell}\} \cup Z$.  If $x_{1}\neq y_{1}$, then $\pi_{x_{1},y_{1}}\in \fip$ since $\fip$ is left-compressed.  Then $\fix(\pi_{x_{1},y_{1}})=\{x_{1},y_{2},\dots,y_{\ell}\} \cup \{z \in Z : z \neq y_{1}\}$.  Repeating this in order of increasing $i$ for all pairs $(x_{i},y_{i})$ gives a permutation $\pi\pr \in \fip$ such that $\fix(\pi\pr)=\{x_{1},x_{2},\dots,x_{\ell}\} \cup \{z \in Z : z \neq y_{i} \mbox{ for any } i \in [\ell]\}$.  It follows that $\fix (\pi\pr) \subseteq \fix(\sigma)$.  Since $\pi\pr \in \fip$, there exists some set $B\pr\pr \in g(\fip)$ such that $\pi\pr \in \mathscr U_{p}(B\pr\pr)$.  It then follows that $B\pr\pr \subseteq \fix(\pi\pr) \subseteq \fix(\sigma)$ and that $\sigma \in \mathscr U_{p}(B\pr\pr)$.  Hence $\sigma \in \fip$ as required.
\qed

\begin{lemma}\label{t+1}
Let $\fip \in I(n,t)$ be maximal, compressed and fixed, and let $n> t+1$.  Let $E_{1}$ and $E_{2}$ be sets in $g(\fip) \in G_{*}(\fip)$.  If there exist some $i<j \in [n]$ such that $i \notin E_{1} \cup E_{2}$ and $j \in E_{1} \cap E_{2}$, then $\left| E_{1} \cap E_{2} \right| \geq t+1$.
\end{lemma}

\proof
Suppose there is such an $i,j$ pair.  Then either $(E_{1}\backslash\{j\}) \cup \{i\}$ or a proper subset thereof will be in $g(\fip)$.  Call this set $B$.  If $i \notin B$, then $B \subset E_{1}$ which is not possible since $g(\fip)$ is minimal by inclusion.  Therefore we conclude that $i \in B$.  It follows from Lemma~\ref{gen set int} that $B$ and $E_{2}$ must $t$-intersect.  Since $j \notin B$ and $i \notin E_{2}$, the result follows.
\qed

\begin{lemma}\label{disjoint union}
Let $\fip \in I(n,t)$ be maximal, fixed and compressed and let $g(\fip)$ be a generating set in $G_{*}(\fip)$.  For a set $E \in g(\fip)$, define $$\mathscr D(E)=\{\sigma \in \sym : \fix(\sigma) \cap [s^{+}(E)]= E\} \: .$$
Then $\fip$ is a disjoint union $$\fip= \dot{\bigcup_{E\in g(\fip)}} \mathscr D(E) \: .$$ 
\end{lemma}  

\proof
We first show that every permutation in $\fip$ is contained in $\mathscr D(E)$ for some $E \in g(\fip)$ and then show by contradiction that $\mathscr D(E_{1}) \cap \mathscr D(E_{2}) = \emptyset$ for any $E_{1}\neq E_{2} \in g(\fip)$.

Let $\sigma\! \in\! \fip$.  Then there exists some $E\! \in\! g(\fip)$ such that $\sigma\! \in\! \up(E)$.  If $\fix(\sigma) \cap [s^{+}(E)]= E$, then $\sigma \in \mathscr D(E)$, so assume that $\fix(\sigma) \cap [s^{+}(E)]\neq E$.   Let $\fix(\sigma)\cap [s^{+}(E)] = \{b_{1},b_{2},\dots,b_{m}\}$ where $b_{1}\!<\! b_{2}\!<\!\dots\!<\!b_{m}\!=\!s^{+}(E)$, and let $\left| E\right| = k<m$.  Then $E \subsetneq \{b_{1},b_{2},\dots,b_{m}\}$ and $B=\{b_{1},b_{2},\dots,b_{k}\} \in \mathscr L(E)$.  By Lemma~\ref{properties}, there is some $B\pr \subseteq B$ such that $B\pr \in \mathcal L_{*}(g(\fip))=g(\fip)$.  By construction $\fix(\sigma) \cap [s^{+}(B\pr)]=B\pr$ and so $\sigma \in \mathscr D(B\pr)$.

Suppose there is some $\pi \in \fip$ such that $\pi \in \mathscr D(E_{1}) \cap \mathscr D(E_{2}) \neq \emptyset$ for some $E_{1}\neq E_{2} \in g(\fip)$. Then $\fix(\pi) \cap [s^{+}(E_{1})]=E_{1}$  and $ \fix(\pi) \cap [s^{+}(E_{2})]=E_{2}$.
Clearly, $s^{+}(E_{1}) \neq s^{+}(E_{2})$ since $E_{1} \neq E_{2}$, so assume without loss of generality that $s^{+}(E_{1}) > s^{+}(E_{2})$.  Then $E_{2} \subsetneq E_{1}$ and hence $g(\fip$) is not minimal with respect to set inclusion which is a contradiction of $g(\fip) \in G_{*}(\fip)$.
\qed

\begin{lemma}\label{size of D}
Let $\fip \in I(n,t)$ be maximal, fixed and compressed and let $g(\fip) \in G_{*}(\fip)$.  Choose a set $\hat{E} \in g(\fip)$ such that $s^{+}(\hat{E})=s^{+}(g(\fip))$.  Then $\mathscr D(\hat{E})$ is the set of all permutations in $\fip$ which are generated by $\hat{E}$ alone.  That is, $$\mathscr D(\hat{E}) =\mathscr U_{p}(\hat{E})\,\backslash\, \mathscr U_{p}\,(g(\fip)\backslash\{\hat{E}\}) \: .$$
Further, $$ \left | \mathscr D(\hat{E}) \right | =  \sum_{j=0}^{s^{+}(\hat{E})-|\hat{E}|}(-1)^{j}\binom{s^{+}(\hat{E})\!-\!|\hat{E}|}{j}\left(n\!-\!|\hat{E}|\!-\!j\right)! \: .$$ 
\end{lemma}
\proof
Since $\mathscr D(E) \subseteq \up(E)$ for all $E \in g(\fip)$, Lemma~\ref{disjoint union} implies that any permutation in $\fip$ generated only by $\hat{E}$ will be in $\mathscr D(\hat{E})$.  Since $g(\fip)$ is minimal with respect to set inclusion, any $E \in g(\fip)\backslash\{\hat{E}\}$ will contain some $x \in [n]$ such that $x \notin \hat{E}$ and $x < s^{+}(\hat{E})$.  Hence $\mathscr D(\hat{E})$ will not contain any permutations generated by $g(\fip)\backslash\{\hat{E}\}$.

The formula for the cardinality of $\mathscr D(\hat{E})$ is derived by using the principle of inclusion-exclusion to count the number of permutations in $\sym$ that have all the elements of $\hat{E}$ fixed and none of the elements of $[s^{+}(\hat{E})]\backslash \hat{E}$ fixed.  
\qed 

\begin{lemma}\label{D prime}
Let $\fip \in I(n,t)$ be maximal, fixed and compressed and let $g(\fip) \in G_{*}(\fip)$. Let $\hat{E} \in g(\fip)$ such that $s^{+}(\hat{E})=s^{+}(g(\fip))$ and let ${\hat{E}}\pr = \hat{E}\backslash\{s^{+}(\hat{E})\}$.  Define $$\mathscr D\pr(\hat{E})= \{\sigma \in \sym : \fix(\sigma) \cap [s^{+}(\hat{E})\!-\!1]= {\hat{E}}\pr\}\:.$$  Then $\mathscr D\pr(\hat{E})$ will be the set of all permutations which are generated by ${\hat{E}}\pr$ and not by $g(\fip)\backslash \{\hat{E}\}$.  Furthermore, $$\left| \mathscr D\pr(\hat{E}) \right| \geq \left (n\!-\!|\hat{E}| +1\right) \left|\mathscr D(\hat{E})\right|\: .$$ 
\end{lemma}
\proof
From the definitions of $\mathscr D(\hat{E})$ and $\mathscr D\pr(\hat{E})$, it is easy to see that $$\mathscr D\pr(\hat{E}) = \mathscr D(\hat{E})\: \dot{\cup} \: \{\sigma \in \sym : \fix(\sigma) \cap [s^{+}(\hat{E})]=\hat{E}\pr\}\:,$$  that is, we can partition $\mathscr D\pr(\hat{E})$ based on whether or not the permutations fix $s^{+}(\hat{E})$.  We will show that neither partition contains permutations generated by $g(\fip)\backslash \{\hat E\}$.

Lemma~\ref{size of D} gives the desired result for $\mathscr D(\hat{E})$.  To show that $\{\sigma \in \sym : \fix(\sigma) \cap [s^{+}(\hat{E})]=\hat{E}\pr\}$ does not contain any permutations generated by $g(\fip)\backslash\{\hat{E}\}$, recall from the proof of Lemma~\ref{size of D} that any $E \in g(\fip)\backslash\{\hat{E}\}$ will contain some $x \in [n]$ such that $x \notin \hat{E}$ and $x < s^{+}(\hat{E})$.  Clearly $x \notin \hat{E}$ implies that $x \notin {\hat{E}}\pr$ and thus $\{\sigma \in \sym : \fix(\sigma) \cap [s^{+}(\hat{E})]=\hat{E}\pr\}$ does not contain any permutations generated by $g(\fip)\backslash\{\hat{E}\}$.

It remains to be shown that all of the permutations generated by ${\hat{E}}\pr$ and not by $g(\fip)\backslash\{\hat{E}\}$ are in $\mathscr D\pr(\hat{E})$.  Suppose there is a permutation, $\sigma$, such that $\sigma\! \in\! \mathscr U_{p}({\hat{E}}\pr)\,\backslash\,\mathscr U_{p}\,(g(\fip)\backslash\{\hat{E}\})$ and $\sigma \notin \mathscr D\pr(\hat{E})$.  Then ${\hat{E}}\pr \subseteq \fix(\sigma)$ and $(\fix(\sigma) \cap [s^{+}(\hat{E})\!-\!1]) \neq {\hat{E}}\pr$.  Let $\fix(\sigma) \cap [s^{+}(\hat{E})\!-\!1] = \{b_{1},b_{2},\dots,b_{\ell}\}$ where $b_{1}< b_{2}<\dots<b_{\ell}$ and let $|\hat{E}| = k$.  Then $k \leq \ell$ since $\fix(\sigma)$ must contain some $x \in [s^{+}(\hat{E})\!-\!1]$ such that $x \notin {\hat{E}}\pr$.  Since $g(\fip) \in G_{*}(\fip)$, there is some $A \subseteq \{b_{1},b_{2},\dots,b_{k}\}\subseteq \fix(\sigma)$ such that $A \in g(\fip)$.  This implies that $\sigma \in \up(A)$ which contradicts the assumption that $\sigma\! \in\! \mathscr U_{p}({\hat{E}}\pr)\,\backslash\,\mathscr U_{p}\,(g(\fip)\backslash\{\hat{E}\})$.   

Let $\sigma$ be any permutation in $\mathscr D(\hat{E})$.  Then $\sigma \in \mathscr D\pr(\hat{E})$.  For any $x \in [n]\backslash\hat{E}$, the permutation formed by transposing $\sigma(x)$ and $s^{+}(\hat{E})$ will be an element of $\mathscr D\pr(\hat{E})\backslash\mathscr D(\hat{E})$.  Therefore, $$\left| \mathscr D\pr(\hat{E})\right| \geq \left|\mathscr D(\hat{E})\right| +\left (n\!-\!|\hat{E}|\right) \left|\mathscr D(\hat{E})\right| =\left (n\!-\!|\hat{E}| +1\right)\left|\mathscr D(\hat{E})\right| \: .\qed$$

\begin{corollary}\label{strictineq}
Let $\fip \in I(n,t)$ be maximal, fixed and compressed and let $g(\fip)$ be a generating set in $G_{*}(\fip)$. Let $\hat{E}$ be a set in $g(\fip)$ such that $s^{+}(\hat{E})=s^{+}(g(\fip))$ and let ${\hat{E}}\pr = \hat{E}\backslash\{s^{+}(\hat{E})\}$ and define $\mathscr D\pr(\hat{E})$ as in Lemma~\ref{D prime}. If $s^{+}(\hat{E})\!-\! | \hat{E} | \geq 1$, then 
\begin{equation}\label{eq:ineq}
\left| \mathscr D\pr(\hat{E}) \right| > \left(n\!-\! |\hat{E}| +1\right) \left|\mathscr D(\hat{E})\right|\:.
\end{equation}
\end{corollary}
\proof
We will show that the inequality is strict when $s^{+}(\hat{E})\!-\!| \hat{E} | \geq 1$ by identifying a permutation in $\mathscr D\pr(\hat{E})$ that is not formed by transposing two elements of a permutation from $\mathscr D(\hat{E})$ and so is not counted in Lemma~\ref{D prime}.

If $s^{+}(\hat{E})\!-\!|\hat{E}| \geq 1$, there will be some $y \in [s^{+}(\hat{E})]$ such that $y \notin \hat{E}$.  Let $\pi$ be a permutation from $\sym$ such that $\fix(\pi)=\hat{E}\cup\{y\}$.  Then $\pi \notin \mathscr D(\hat{E})$.  Now consider the permutation $\tilde{\pi}=(y,s^{+}(\hat{E}))\comp \pi$.  Then $\tilde{\pi}(y)=s^{+}(\hat{E})$ and $\tilde{\pi}(s^{+}(\hat{E}))=y$.  Thus $\fix(\tilde{\pi})={\hat{E}}\pr$ and so $\tilde{\pi} \in \mathscr D\pr(\hat{E})$ but $\tilde{\pi} \notin \mathscr D(\hat{E})$.  Equation~\ref{eq:ineq} then follows. 
\qed

\section{Proof of Theorem~\ref{KR with bound}}


We now use the results established in Section 2.3 to prove the following lemma.
\begin{lemma}\label{main lemma}

Let $\fip \in I(n,t)$ be a compressed fixed family of maximum possible size.  If $n \geq 2t + 1$, then $s_{min}(G(\fip)) \leq t$.

\end{lemma}

\proof
By Lemma~\ref{not empty} we have that $G(\fip) \neq \emptyset$ and by Lemma~\ref{properties} that $G_{*}(\fip) \neq \emptyset$ and$s_{min}(G(\fip)) = s^{+}(g(\fip))$ for some $g(\fip) \in G_{*}(\fip))$.  We will use this $g(\fip)$ throughout the proof.   

Assume that $s^{+}(g(\fip)) = t + \delta$ for some positive integer $\delta$.  Clearly $t+\delta \leq n$.  Now partition $g(\fip)$ into two disjoint collections of sets,
 $$g_{0}(\fip)= \left\{B \in g(\fip) : s^{+}(B)= t + \delta\right\}\:,$$ \vspace{-12pt}and 
 $$g_{1}(\fip)= g(\fip)\backslash \, g_{0}(\fip)\:.$$
For any $B\! \in \! g_{0}(\fip)$ and $A\! \in \! g_{1}(\fip)$, it is clear that $\left| (B\backslash\{t\!+\!\delta\}) \cap A \right | \! \geq t \!\; .$
  
Next, we partition the sets in $g_{0}(\fip)$ according to their cardinality.  Let $$\rpart{i}=\left\{B \in g_{0}(\fip) : \left|B\right| = i\right\}\: .$$

If $i \leq t$, then $\rpart{i}= \emptyset$.  To see why this is true, assume there is a set $B \in \rpart{i}$ for some $i \leq t$.  Then $B$ will be a set containing $t+\delta$ with cardinality less than or equal to $t$.  Since $g(\fip)$ is left shifted, it follows that $[i] \in g(\fip)$.  But $\left| [i] \cap B \right| < t$ and so it follows from Lemma~\ref{gen set int} that $\fip$ is not $t$-cycle-intersecting.  Also, if $i=t + \delta$, then $g(\fip)=\{[t+\delta]\}$ and $\fip$ is clearly not a $t$-cycle-intersecting family of maximum size.  Hence,
$$g_{0}(\fip) = \dot{\bigcup_{t<i<t+\delta}} \mathscr R_{i} \: .$$
Also, if $\rpart{i} \neq \emptyset$, then $i \leq n-2$ since by definition there are no sets of cardinality $n\!-\!1$ in $g(\fip)$ and $i \neq n$ is implied by $i < t+\delta.$

We now consider the set system formed by removing $t+\delta$ from the sets in $g_{0}(\fip)$.  Let $$\rprime{i} = \left\{E\backslash\{t+\delta\} : E \in \rpart{i}\right\}\:.$$
Clearly $\left| \rpart{i}\right| = \left | \rprime{i}\right|$ and $\left|E\pr\right|=i-1$ for $E\pr = E\backslash\{t+\delta\} \in \rprime{i}$.

For any $E_{1}, E_{2}\in g_{0}(\fip)$ with $\left|E_{1} \cap E_{2}\right|=t$, it follows from Lemma~\ref{t+1} that $i \in E_{1} \cup E_{2}$ for all $i < t+\delta$.  A simple counting argument then gives us $$\left|E_{1}\right| + \left| E_{2}\right| = \left| E_{1} \cup E_{2}\right| + \left| E_{1} \cap E_{2}\right| =2t + \delta \: .$$
Therefore, if $\left|E_{1}\right| + \left| E_{2}\right| \neq 2t + \delta$, then $\left|E_{1} \cap E_{2}\right| > t$.  It then follows that for any $E\pr_{1} \in \rprime{i}$ and $E \pr_{2} \in \rprime{j}$ with $i+j \neq 2t + \delta$, $$\left|E\pr_{1} \cap E\pr_{2}\right| \geq t \: .$$

We claim that for $i,j \in \mathbb N$ such that $i + j =2t +\delta$, if $\rpart{i} \neq \emptyset$, then $\rpart{j} \neq \emptyset$.\\
To prove this claim, suppose that there is a set $E_{1}$ in $\rpart{i}$ such that $\left|E_{1} \cap E\right|\geq t+1$ for all $E \in g_{0}(\fip)$.  Then $\left|(E_{1}\backslash\{t+\delta\}) \cap E\right| \geq t$ for all $E \in g_{0}(\fip)$.  Since $E_{1}\backslash\{t+\delta\}$ also $t$-intersects with the sets in $g_{1}(\fip)$ and $\fip$ is maximal, we must have  
$E_{1}\backslash\{t+\delta\} \in g(\fip)$.  But this implies that $E_{1} \notin g(\fip)$ since $g(\fip)$ is minimal by inclusion, a contradiction.  Therefore, for any $E_{1} \in \rpart{i}$, there must exist at least one other set, $E_{2} \in g_{0}(\fip)$, such that $\left| E_{1} \cap E_{2} \right| = t$.  It then follows that $\left|E_{1}\right| + \left| E_{2}\right|= 2t+\delta$.  Thus, if $\rpart{i} \neq \emptyset$, then $\rpart{j} \neq \emptyset$ where $j= 2t + \delta - i$.

We have already seen that $\rpart{i}= \emptyset$ when $i\leq t$ or $i \geq t+\delta$.  We now show that if $\fip$ is as large as possible and $n \geq 2t+1$, then $\rpart{i} = \emptyset$ for all $i$ and hence $s^{+}(g(\fip))=t$.  We will do this by assuming that $\rpart{i} \neq \emptyset$ for some $i \in \{t+1, \dots , t +\delta -1\}$ and then constructing a family of permutations that is larger than $\fip$.  We will consider the cases where $i \neq 2t+\delta -i$ and where $i = 2t+\delta -i$ separately.

\begin{enumerate}
\item  
Case 1:  $i \neq 2t+\delta -i$.\\
Recall that $\rpart{i} \neq \emptyset$ implies that $\rpart{2t + \delta -i} \neq \emptyset$ and that $i \leq n-2$ and $2t+\delta-i \leq n-2$.  Consider the sets $$f_{1}=\left(g(\fip)\backslash(\rpart{i}\cup\rpart{2t+\delta-i})\right)\cup \rprime{i}$$ and $$f_{2}=\left(g(\fip) \backslash(\rpart{i}\cup\rpart{2t+\delta-i})\right)\cup\rprime{2t+\delta-i} \: .$$
We first show that $f_{1}$ is $t$-intersecting.  Clearly, $g(\fip) \backslash(\rpart{i}\cup\rpart{2t+\delta-i})$ is a $t$-intersecting set system since by Lemma~\ref{gen set int}, $g(\fip)$ is $t$-intersecting.  Let $E\pr_{1}$ be a set in $\rprime{i}$.  Then $E_{1}=E\pr_{1}\cup \{t+\delta\}$ is a set in $\rpart{i}$.  As shown previously, $\left|E_{1} \cap E_{j}\right| \geq t+1$ for all $E_{j} \in \rpart{j}$ where $j \neq 2t+\delta-i$.  Hence, $E\pr_{1}$ will $t$-intersect with all sets in $g_{1}(\fip)$, $g_{0}(\fip)\backslash \rpart{2t+\delta-i}$ and $\rprime{i}$.  Thus $f_{1}$ is a $t$-intersecting set system.

A similar argument can be used to show that $f_{2}$ is $t$-intersecting.

Let $\mathcal B_{1}= \up(f_{1})$ and let $\mathcal B_{2}=\up(f_{2})$.  Then $\mathcal B_{1}$ and $\mathcal B_{2}$ are $t$-cycle-intersecting families of permutations.  We claim that $$\max_{i=1,2} \left |\mathcal B_{i}\right| > \left|\fip\right| \: .$$

For any set $E \in \rpart{i} \subseteq g_{0}(\fip)$, consider $\mathscr D(E)$ as defined in Lemma~\ref{disjoint union}.  Recall that $\mathscr D(E_{1})\cap\mathscr D(E_{2}) = \emptyset$ for all $E_{1},E_{2} \in g(\fip)$.  By Lemma~\ref{size of D} we have that $\mathscr D(E)$ is the set of all permutations generated only by $E$.  Therefore, the permutations generated by $\rpart{i}$ and not by $g(\fip)\backslash \rpart{i}$ will be given by $$\mathscr D(\rpart{i}) = \dot{\bigcup_{E\in\rpart{i}}} \mathscr D(E)  \: .$$
Again by Lemma~\ref{size of D}, $$ \left | \mathscr D(E) \right | =  \sum_{j=0}^{s^{+}(E)-|E|}\left(-1\right)^{j}\binom{s^{+}(E)-|E|}{j}\left(n-|E|-j\right)! \: .$$
Since $\left|E\right| = i $ and $s^{+}(E) = t +\delta$ for all $E\in \rpart{i}$, $$\left|\mathscr D(\rpart{i})\right| = \left|\rpart{i}\right|\cdot\left|\mathscr D(E)\right|= \left|\rpart{i}\right|\sum_{j=0}^{t+\delta-i}\left(-1\right)^{j}\binom{t+\delta-i}{j}\left(n-i-j\right)!\:.$$
Similarly, $$\mathscr D(\rpart{t+\delta-i}) =\dot{\bigcup_{E\in\rpart{t+\delta-i}}} \mathscr D(E)$$ is the set of permutations generated by $\mathscr D(\rpart{t+\delta-i})$ and not by $g(\fip)\backslash\mathscr D(\rpart{t+\delta-i})$.
 
It follows from Lemma ~\ref{D prime} that $$\mathscr D\pr(\rpart{i})=\dot{\bigcup_{E\in\rpart{i}}} \mathscr D\pr(E)$$ is the set of all permutations generated by $\rprime{i}$ and not by $g(\fip)\backslash\rpart{i}$.
Since $i < t+\delta$, we can use Equation~\ref{eq:ineq} from Corollary~\ref{strictineq} to show that
\begin{align*}
\left|\mathcal B_{1}\right| &= \left|\fip\right|-\left(\left|\mathscr D(\rpart{i})\right|+ \left|\mathscr D(\rpart{2t+\delta-i})\right|\right) + \left| \mathscr D\pr(\rpart{i})\right|\\
&>  \left|\fip\right| - \left| \mathscr D(\rpart{2t+\delta-i})\right| + \left(n-i\right)\cdot\left|\mathscr D(\rpart{i})\right|\\
\intertext{and}
\left|\mathcal B_{2}\right| &= \left|\fip\right|-\left(\left|\mathscr D(\rpart{i})\right|+ \left|\mathscr D(\rpart{2t+\delta-i})\right|\right) + \left| \mathscr D\pr(\rpart{2t+\delta-i})\right|\\
&> \left|\fip\right| - \left| \mathscr D(\rpart{i})\right| + \left(n\!-\!2t\!-\!\delta +i\right)\cdot\left|\mathscr D(\rpart{2t+\delta-i})\right|\:.
\end{align*}

We now prove by contradiction that either $\mathcal B_{1}$ or $\mathcal B_{2}$ is larger than $\fip$.\\  If $\left| \fip \right| \geq \left| \mathcal B_{1}\right|$, then 
$$\left|\mathscr D(\rpart{2t+\delta-i})\right| > (n-i)\cdot\left|\mathscr D(\rpart{i})\right|$$ and if $\left| \fip \right| \geq \left| \mathcal B_{2}\right|$, then $$\left|\mathscr D(\rpart{i})\right| > (n-2t-\delta+i)\cdot\left|\mathscr D(\rpart{2t+\delta-i})\right|\:.$$
Thus, if $\fip$ is larger than both $\mathcal B_{1}$ and $\mathcal B_{2}$, $$\left| \mathscr D(\rpart{2t+\delta-i})\right| > (n-i)(n\!-\!2t\!-\!\delta+i)\cdot\left|\mathscr D(\rpart{2t+\delta-i})\right|$$ and $$ 1 > (n\!-\!2t\!-\!\delta+i)(n-i)\:.$$


However, both $n-i \geq 2$ and $n-(2t+\delta-i) \geq 2$.  Therefore at least one of $\mathcal B_{1}$ or $\mathcal B_{2}$ will be larger than $\fip$.  Thus if $\fip$ is of maximum size, then $\rpart{i}=\rpart{2t+\delta-i}=\emptyset$ for all $i \neq 2t + \delta -i$.

\item
Case 2:  $i= 2t+\delta - i$.\\
In this case, $i= t + \frac{\delta}{2}$ and hence $\delta$ must be divisible by $2$.  Consider $\rprime{t+ \delta /2}$.  For any 
$B\pr \in \rprime{t+\delta/2}$, note that $\left| \, B\pr \right| = t + \frac{\delta}{2} - 1$  and $B\pr \subset [t\!+\!\delta \!-\!1]$.  
First we claim that there exists some $a \in [t\!+\!\delta\!-\!1]$ and some collection of subsets, $\mathscr T\pr \subseteq \rprime{t + \delta/2}$, such that for all $B\pr \in \mathscr T\pr\,$, we have that $a \notin B\pr$ and 
\begin{equation}
\label{eq:subsetsize}
\left| \mathscr T\pr  \right| \geq \left| \rprime{t+\delta/2} \right| \cdot \frac{\delta/2}{t\!+\!\delta\!-\!1}=\left| \rpart{t+\delta/2} \right| \cdot \frac{\delta}{2(t\!+\!\delta\!-\!1)}.
\end{equation}
To prove this claim we will consider the complements of the sets in $\rprime{t+\delta/2}$ in $[t+\delta\!-\!1]$.  For all $B\pr_{j} \in \rprime{t+\delta/2}$, let $C\pr_{j} = [t+\delta \!-\!1]\backslash B\pr_{j}$.  Then $ \left| C\pr_{j} \right| = \frac {\delta}{2}$ and there will be $\left | \rprime{t+\delta/2}\right |$ distinct sets of size $\frac{\delta}{2}$ with entries from $[t+\delta \!-\!1]$.  By the pigeonhole principle, there exists some $a \in [t+\delta \!-\!1]$ such that $a$ is in at least $\frac{| \rprime{t+\delta/2}|}{t+\delta \!-\!1}\cdot\frac{\delta}{2}$ of the $C\pr_{j}\:$s.  The claim then follows easily.

For any $B\pr_{1}, B\pr_{2} \in \mathscr T\pr$, Lemma~\ref{t+1} gives $$ \left | B\pr_{1} \cap B\pr_{2} \right | \geq t.$$  From Case 1 we have that $\rpart{i} = \emptyset$ for $i \neq t + \frac{\delta}{2}$.  Hence $$f\pr = \left(g(\fip)\backslash \rpart{t+\delta/2}\right) \cup \mathscr T\pr$$ is a $t$-intersecting set system and therefore $\up(f\pr)$ will be a $t$-cycle-intersecting family of permutations.

We now show that $\left |\up(f\pr) \right | > \left |\fip\right |$ when $n \geq 2t + 1$.\\  Let
\begin{eqnarray*}
\mathscr D_{1} &=& \up\,\left(g(\fip)\,\backslash \, \rpart{t+\delta/2}\right),\\
\mathscr D_{2} &=& \up\,(\rpart{t+\delta/2})\,\backslash\, \mathscr D_{1}\:,\\
\mathscr D_{3} &=& \up\,(\mathscr T\pr)\,\backslash \, \mathscr D_{1}\:.\\
\end{eqnarray*}
Then $\fip = \mathscr D_{1} \, \dot\cup \: \mathscr D_{2}$ and $\up(f\pr) = \mathscr D_{1}\, \dot\cup \: \mathscr D_{3}$.  Hence to show that $\left |\up(f\pr) \right | > \left |\fip\right |$, it is sufficient to show that $\left | \mathscr D_{3} \right | > \left| \mathscr D_{2} \right |$.

From Lemma~\ref{size of D} we have
\begin{equation}
\label{eq:D2}
\left | \mathscr D_{2} \right | = \left |\rpart{t+\delta/2}\right | \sum_{j=0}^{\delta/2}(-1)^j\binom{\delta/2} { j}(n\!-\!\left(t+\delta/2)\!-\!j\right)!\,.
\end{equation}

Now we determine a lower bound on the size of $\mathscr D_{3}$.  Let $B \in \rpart{t+\delta/2}$.  Recall that $\mathscr D(B)$ is the set of all permutations generated by $B$ and not by $g(\fip)\backslash B$.  As shown in Lemma~\ref{D prime}, for each permutation in $\mathscr D(B) $ there will be at least $n-i+1= n-(t+\delta/2)+1$ permutations generated by $B\pr = B\backslash \{t+\delta\}$ and not by $g(\fip)\backslash B$.  From Lemma~\ref{size of D} we have $$\mathscr D(B) = \sum^{(t+\delta)- (t+\delta/2)}_{j=0}(-1)^{j}\binom{(t + \delta)-(t+ \delta/2)}{j} \left(n\!-\! (t+\delta/2) \!-\!j\right)!\,.$$  As in Case 1, we use Equation~\ref{eq:ineq} from Corollary~\ref{strictineq}, to show that $$\left | \mathscr D_{3} \right | > \left | \; \mathscr T \pr \right | \left(n\!-\!(t+\delta/2)+1\right) \sum^{\delta/2}_{j=0}(-1)^{j}\binom{\delta/2}{j} \left(n\!-\! (t+\delta/2) \!-\!j\right)!\:.$$  Combining this with Equations~\ref{eq:subsetsize} and~\ref{eq:D2} gives

$$\left| \mathscr D_{3}\right| >\frac{\delta}{2(t+\delta+1)}\left(n\!-\!t\!-\!\delta/2+1\right)\cdot \left|\mathscr D_{2}\right|\:.$$  Hence it is sufficient to show that $$\frac{\delta}{2(t+\delta\!-\!1)}\cdot\left(n\!-\!t\!-\!\delta/2+1\right) \geq 1$$ or equivalently that
\begin{equation}
\label{eq:quad}
0 \geq \delta^{2} +\delta(2+2t-2n)+4t-4 \:.
\end{equation}
Since $n \geq t+\delta$ and $\delta \geq 2$,  we have $2 \leq \delta \leq n-t$ .
Thus if the inequality in Equation~\ref{eq:quad} holds for $\delta=2$ and $\delta=n-t$, it will hold for all possible values of $\delta$.

If $\delta=2$, Equation~\ref{eq:quad} is equivalent to $n \geq 2t+1$.

Now let $\delta = n-t$.  Then Equation~\ref{eq:quad} becomes 
\begin{equation}
\label{eq:last}
0 \geq n(-n+2t+2) + (2t-t^{2}-4)\:.
\end{equation}
Clearly $(2t-t^{2}-4) <0$ for all $t \in \mathbb N$ and $n(-n+2t+2) \leq 0$ provided that $n\geq 2t+2$.  Thus it is easy to see that the inequality in  Equation~\ref{eq:last} holds provided that $n \geq 2t+2$. What if $n=2t+1$?  Then $\delta = t+1$ and $t$ must be odd since $\delta$ is even.  Substituting $2t+1$ for $n$ in Equation~\ref{eq:last} gives $$0 \geq 4t-t^{2}-3\,.$$  This inequality holds for $t=1$ and $t \geq 3$ and thus it holds for all odd $t \in \mathbb N$.

Therefore the inequality in Equation~\ref{eq:quad} holds for all $t \in \mathbb N$ when $n \geq 2t+1$.  Hence $\left|\mathscr D_{3}\right| > \left|\mathscr D_{2}\right|$ and thus $\left|\mathscr U_{p}(f\pr)\right| > \left| \fip \right|$ when $n\geq 2t+1$.
\end{enumerate}

We have now shown that  when $i \neq 2t + \delta -i$ and when $i = 2t+ \delta -i$, if $s_{min}(G(\fip))>t$, then $\fip$ does not have the maximum possible size. \qed
 
Lemmas 3.1 and 2.11 are sufficient to prove Theorem~\ref{KR with bound} if the theorem is restricted to fixed, compressed $t$-cycle-intersecting families of permutations.  However, as shown in Sections 2.1 and 2.2, any maximal $t$-cycle-intersecting family of permutations can be transformed into a fixed and compressed $t$-cycle-intersecting family of the same size and the resulting family will be the stabilizer of $t$ points only if the original family was the stabilizer of $t$ points.  This completes the proof of Theorem~\ref{KR with bound}.

\section{Discussion}


We have claimed that $n \geq 2t + 1$ is the best possible lower bound on $n$.  To see that this is true, let $\fip_{0}$ be the pointwise stabilizer of $[t]$ and consider the family of permutations, $\fip_{1}$, consisting of all permutations that fix at least $t+1$ of the integers $1,2,\dots,t+2$.   That is, $$\fip_{1}= \{\sigma \in \sym : \left|\fix \sigma \cap [t+2] \right| \geq t+1\}\,.$$  If $n=2t$, then $\left|\fip_{0}\right|=t!=(t-2)!(t^{2}-t)$ and $\left|\fip_{1}\right|=(t-2)!(t^{2}-3)$.  Clearly, $\left| \fip_{1}\right| =\left| \fip_{0}\right|$ when $t=3$ and $\left|\fip_{1}\right| > \left|\fip_{0}\right|$ when $t >3$.  In fact, $\left|\fip_{1}\right| > \left|\fip_{0}\right|$ for all $n$ such that $t+3 \leq n < 2t+1$ with $t \neq 3$ (see~\cite{alison}). 

Ellis et al~\cite{MR2784326} conjecture that a $t$-intersecting family of permutations of maximum size will have the form $\sigma\fip_{i}\tau$ where $\sigma,\tau \in \sym$ and $$\fip_{i} = \{\sigma \in \sym : \left| \fix\sigma \cap [t+2i]\right| \geq t+i \}$$ for some $0\leq i \leq (n-k)/2$.  In this paper, we have adapted only a portion of the Ahlswede and Khachatrian proof of the Complete \erd-Ko-Rado theorem and it is possible that other techniques in their proof can be used to prove this conjecture for $t$-cycle-intersecting families of permutations.

The main problem encountered in applying the method used in this paper to $t$-intersecting families of permutations is the fixing operation.  Applying the $ij$-fixing operation does not preserve the $t$-intersection.  The $x$-fixing operation introduced in~\cite{MR2009400} does preserve the $t$-intersection but does not preserve the size; the family obtained by applying the $x$-fixing operation may contain fewer permutations than the original family.  A more detailed examination of these operations can be found in~\cite{alison}.  It is interesting to note that a family of permutations is closed under the $ij$-fixing operation if and only if it is closed under the $x$-fixing operation.  Also, families of the form $\sigma \fip_{i}\tau$ as defined above can be transformed into fixed families without loss of $t$-intersection or size.  Thus a proof of the conjecture would allow results for fixed $t$-intersecting families of maximum size to be extended to all $t$-intersecting families of maximum size.

\end{document}